\documentclass{amsart}
\usepackage{amssymb}
\newtheorem{Theo}{Theorem}
\newtheorem{cor}{Corollary}

\newcommand{\N}{\mathbb{N}}
\newcommand{\Q}{\mathbb{Q}}
\begin{document}
\title[Character sums and Bianchi groups]{Sifted character sums  and
  free quotients of Bianchi groups} 
\author[J.-C. Schlage-Puchta]{Jan-Christoph Schlage-Puchta}
\begin{abstract}
We show that the Bianchi group $\mathrm{PSL}_2(\mathcal{O}_d)$, where
$\mathcal{O}_d$ is the ring of integers in $\Q(\sqrt{d})$, $d<0$, has a
free quotient of rank $\geq |d|^{\frac {1}{4}-\epsilon}$, as
$|d|\to\infty$. To do so, we
give an estimate for a sifted character sum.
\end{abstract}
\maketitle
Let $d<0$ be a fundamental discriminant, $\mathcal{O}_d$ be the ring of
integers in $\Q(\sqrt{d})$, and
$\Gamma_d=\mathrm{PSL}_2(\mathcal{O}_d)$ be the corresponding Bianchi
group. Define the Zimmert set $Z_d$ to be the set of all integers $n$
satisfying the following conditions:
\begin{enumerate}
\item $4n^2+3\leq |d|$, and $n\neq 2$;
\item $d$ is a quadratic non-residue modulo $p$ for all odd prime
  factors $p$ of $n$;
\item If $d\not\equiv 5\pmod{8}$, then $n$ is odd.
\end{enumerate}
Denote by $r(d)$ the rank of the largest free quotient of $\Gamma_d$. 
R. Zimmert\cite{Zimmert} proved that $r(d)\geq|Z_d|$. This relation was used by
Mason, Odoni and Stothers\cite{MOS} to show that for $|d|>10^{476}$,
$r(d)\geq 2$, that is, $\Gamma_d$ has a free non-abelian quotient, and
that, as $|d|\to\infty$, we have $r(d)\gg\log d$. The difficulties
in estimating $|Z_d|$ come from the fact that one has to bound sums of
highly imprimitive characters, here, we avoid this problem by
incorporating a sifting device into the character sum. This
approach is similar to the one used in \cite{Sieb}. We will prove the
following general estimate for imprimitive characters. 
\begin{Theo}
\label{thm:main}
Let $\chi\pmod{q}$ be a character, $P$ an integer. Then we have for each
$x$ and parameter $2\leq R\leq x$ the estimate
\[
\sum_{_{n\leq x}\atop{(n, P)=1}} \chi(n) \ll x^{1-\frac{1}{r}}
R^{\frac{1}{r}} q^{\frac{r+1}{4r^2}+\epsilon} +
x^{1+\epsilon}\sum_{{R<t\leq x}\atop{t|P}}\frac{1}{t},
\]
where $r\in\N$ can be chosen arbitrarily, if $q$ is cubefree up to a
factor 8, and $r\in\{1, 2, 3\}$ in general.
\end{Theo}
\begin{cor}
\label{cor:Zimmert}
Let $d<0$ be a fundamental discriminant, $\mathcal{O}_d$ be the ring of
integers in $\Q(\sqrt{d})$. 
Let $\Gamma_d=\mathrm{PSL}_2(\mathcal{O}_d)$ the Bianchi group
corresponding to this ring. Then for each $\epsilon>0$ there is a $d_0$
such that for $|d|>d_0$ the group $\Gamma_d$ maps onto a free group of
rank $\geq |d|^{\frac{1}{4}-\epsilon}$.
\end{cor}
We will make use of the following version of Burgess' famous estimates
for character sums.
\begin{Theo}
\label{thm:Burgess}
Let $\chi\pmod{q}$ be a non-principal character. Then we have for each
$x\geq 2$ the estimate
\[
\sum_{n\leq x}\chi(n) \ll x^{1-\frac{1}{r}} q^{\frac{r+1}{4r^2}+\epsilon},
\]
where $r\in\N$ can be chosen arbitrarily, if $q$ is cubefree, and $r\in\{1, 2, 3\}$
in general.
\end{Theo}
The restriction ``$q$ cubefree'' can be somewhat lessened. If
$q=q_1q_2$ with $q_2$ cubefree, and $(q_1, q_2)=1$, we can write each
character $\chi\pmod q$ as $\chi_1\chi_2$ where $\chi_i$ is a character modulo
$q_i$. Then we have
\[
\sum_{n\leq x} \chi(n)=\sum_{(a, q)=1}\chi_1(a)
\sum_{{n\leq x}\atop{n\equiv a\,\mathrm{mod}\,q_1}}\chi_2(n)\ll q_1 (x/q_1)^{1-\frac{1}{r}} q_2^{\frac{r+1}{4r^2}+\epsilon}.
\]
In particular, if $q$ is cubefree up to a factor 8, we obtain up to a
constant the same estimates as for cubefree integers.

To prove Theorem~\ref{thm:main}, fix a parameter $R$, and decompose
the sum to be estimated as 
\begin{eqnarray*}
\sum_{{n\leq x}\atop{(n, P)=1}} \chi(n) & = & \sum_{n\leq x}\chi(n)
\Big(\sum_{t|(n, P)}\mu(t)\Big)\\
& = & \sum_{n\leq x}\chi(n)
\Big(\sum_{{t\leq R}\atop{t|(n, P)}}\mu(t)\Big) -
\sum_{{n\leq x}\atop{(n, P)>1}}\chi(n)
\Big(\sum_{{t\leq R}\atop{t|(n, P)}}\mu(t)\Big)\\
 & = & \sum\nolimits_1 - \sum\nolimits_2,
\end{eqnarray*}
say. To estimate $\sum_1$, we interchange the order of summation and
obtain
\[
\sum\nolimits_1 \leq \sum_{t\leq R} \left|\sum_{n\leq x/t}
  \chi(n)\right|
\ll \sum_{t\leq R} (x/t)^{1-\frac{1}{r}}
q^{\frac{r+1}{4r^2}+\epsilon}\ll x^{1-\frac{1}{r}}
R^{\frac{1}{r}} q^{\frac{r+1}{4r^2}+\epsilon},
\]
for all $r$ admissible by Theorem~\ref{thm:Burgess}. To bound $\sum_2$,
note that the inner sum vanishes for $2\leq (n, P)\leq R$, hence, we may
restrict the outer sum to terms with $(n, P)>R$. Trivially,
the inner sum is bounded above by $\tau(n)\ll x^\epsilon$, and we obtain
\[
\sum\nolimits_2 \ll x^\epsilon\#\{n\leq x: (n, P)>R\} \leq x^\epsilon
\sum_{{t>R}\atop{t|P}}\left[\frac{x}{t}\right].
\]
The estimates for $\sum_1$ and $\sum_2$ now imply our theorem.

To deduce Corollary~\ref{cor:Zimmert}, fix $c<c'<\frac{1}{4}$,
$0<\epsilon< c'-c$, and let $d$ be a
sufficiently large integer. Since $\Q(\sqrt{-m^2d})=\Q(\sqrt{-d})$, we
may suppose that $d$ is squarefree. Denote by $P$ be the product of all
prime numbers $p$ in the Zimmert set of $d$. Let $\chi$ be the character
defined by $\chi(2)=0$, $\chi(n)=\left(\frac{d}{n}\right)$ for $n$
odd. Then $\chi$ is a character modulo $q=4d$, and $q$ is cubefree up to
a possible factor 8. Setting $R=|d|^c$,
$x=\frac{1}{2}\sqrt{|d|-3}$, and
$r=\left\lceil\frac{1}{1-4c'}\right\rceil$ in Theorem~\ref{thm:main}, we obtain
\begin{eqnarray*}
\sum_{{n\leq x}\atop{(n, P)=1}} \chi(n) & \ll &
|d|^{\frac{1}{2}-(\frac{1}{4}-c')^2} +
x^{1+\epsilon}\sum_{{R<t\leq x}\atop{t|P}}\frac{1}{t}\\
 & \ll & x^{1-\delta} + x^{1-2c'+\epsilon}\#\{R<t\leq x: t|P\}\\
 & \leq & x^{1-\delta} + x^{1-c-c'} |Z_d|.
\end{eqnarray*}
On the other hand, if $n\leq x$ and $(n, Pd)=1$, then
$\chi(n)=1$. Restricting the sum to prime values $n$ we deduce
\[
\sum_{{n\leq x}\atop{(n, P)=1}} \chi(n) \geq 
\sum_{{p\leq x}\atop{{p\nmid P}\atop{p\,\mathrm{prim}}}} \chi(p) \geq \pi(x) - |Z_d|-\omega(|d|),
\]
and comparing these estimates we obtain $|Z_d|\gg |d|^{c}$.

Jan-Christoph Schlage-Puchta\\
Mathematisches Institut\\
Eckerstr. 1\\
79104 Freiburg\\
Germany
\end{document}